\documentclass[12pt]{article}
\usepackage[centertags]{amsmath}
\usepackage{amsfonts}
\usepackage{amssymb}
\usepackage{amsthm}
\usepackage{latexsym}
\usepackage{amsmath}
\usepackage{amsfonts}
\usepackage[active]{srcltx}
\usepackage{comment}
\usepackage{indentfirst}

\usepackage{graphicx}

\DeclareGraphicsExtensions{.eps,.pdf,.jpg,.png}
\usepackage{epstopdf}
\usepackage{subfig}
\usepackage{xcolor}
\usepackage{ulem}
\usepackage{multirow}
\usepackage{adjustbox}

\numberwithin{equation}{section}

\newtheorem{example}{Example}[section]
\newtheorem{theorem}{Theorem}[section]
\newtheorem{nono-lemma}{Lemma}[]
\newtheorem*{lemma}{Lemma}

\begin{document}
\title{Inverse spectral problem for glassy state relaxation approximated by Prony series}

\author{Shuli Chen\thanks{Corresponding author\,%School of Mathematics, Nanjing 210096, P. R. China 
(Email: sli\_chen@126.com).}
\and
Maarten V. de Hoop \thanks{Simons Chair in Computational and Applied
Mathematics and Earth Science, Rice University, Houston TX, USA (Email: mdehoop@rice.edu).}
\and
Youjun Deng\thanks{School of Mathematics and Statistics, Central South University,Changsha 410083, P.R. China (youjundeng@csu.edu.cn).}
\and
Ching-Lung Lin\thanks{Department of Mathematics, National Cheng-Kung University, Tainan 701, Taiwan (Email:
cllin2@mail.ncku.edu.tw).}
\and Gen Nakamura\thanks{Department of
Mathematics, Hokkaido University, Sapporo 060-0808, Japan and Research Center of Mathematics for Social Creativity, Research Institute for Electronic Science, Hokkaido University, Sapporo 060-0812, Japan (Email: gnaka@math.sci.hokudai.ac.jp).}}

\date{}
\maketitle

\begin{abstract}
\noindent The stretched exponential relaxation function is used to analyze the relaxation of the glassy state data. Due to the singularity of this function at the origin, this function is inconvenient for data analysis. Concerning this, a Prony series approximation of the stretched exponential relaxation function (\cite{MM}), which is the extended Burgers model (abbreviated by EBM) known for viscoelasticity equations, was introduced.
In our previous paper \cite{CDDLN}, we gave an inversion method to identify the relaxation tensor of the EBM using clustered eigenvalues of the quasi-static EBM. As a next important research subject of this study, we numerically examine the performance of the inversion method. The performance reveals that it is a powerful method of data analysis, analyzing the relaxation of the glassy state data.


\end{abstract}

\renewcommand{\theequation}{\thesection.\arabic{equation}}

\section{Introduction}\label{intro}
Glassy states can be formed in a variety of materials with different bonding properties such as oxides, alloys, molecules, and polymers. Due to the nonequilibrium of this state, the atoms in these materials exhibit collective motion and exhibit viscoelastic properties upon deformation. For understanding the properties of glasses, it is important to understand the behavior of this deformation (see, for example, [12]).
Non-exponential relaxation of stress is commonly observed for homogeneous glassy materials. More precisely, the following empirical non-exponential function called the stretched exponential function, $G(t)$ given as
\begin{equation}\label{s_ exp_func}
G(t)=\text{exp}\left\{-\left(\frac{t}{\tau}\right)^\beta\right\},  
\end{equation}
is used, where $t\ge 0$ is time, $\tau>0$ is the relaxation time and $0<\beta<1$ is the stretching exponent. Even though this function has a very simple form, it has a singularity at the initial time $t=0$. Since data near the initial time are not accurate, it is a very natural idea to replace this function with a function that approximates this function away from t=0. In fact, the authors of \cite{MM} proposed using a Prony series well-known in viscoelasticity, and demonstrated its effectiveness. Hence, from now on, we replace the stretched exponential function by the following Prony series
\begin{equation}\label{Prony}
G(t)=\sum_{i=1}^N s_i e^{-r_i t},
\end{equation}    
where $r_i>0,\,s_i>0$ for $1\le i\le N$ with $N\in{\mathbb N}$, referred to as {\it dimension}. Here, ${\mathbb N}$ is the set of natural numbers, and we have abused the notation $G(t)$ to denote the Prony series. Also, without loss of generality, we assume that
\begin{equation}\label{I1.3}
0<r_1<r_2<\cdots<r_N
\end{equation}

To describe our inverse spectral problem, we first consider the initial boundary value problem associated with the viscoacoustic wave equation, which has the above mentioned relaxation function, in space dimension one, given as
\begin{equation}\label{IBP}
\left\{
\begin{array}{rcl}
\partial_t^2 u(t, x) &=& \displaystyle{\int_0^t G(t-\tau)\partial_\tau  \partial_x^2 u(\tau, x)\,d\tau}\ \text{in $(0,\infty)\times\Omega$}, \\
u(t,0)&=&0,\quad \partial_x u(t,\frac{\pi}{2})=0\ \text{in $(0,\infty)$},\\[0.25cm]
u(0,x)&=&0\,\,\text{in $\Omega$},
\end{array}
\right.
\end{equation}
with $\Omega = (0,\pi/2)$. Here,
% $u$ denotes the displacement of the beam, and 
we have assumed that the density of the material is constant, equal to $1$, and that $G$ is independent of $x$ for simplicity. Integrating the right-hand side of the equation of motion in \eqref{IBP} by parts, we have
\begin{equation}\label{motion1}
\partial_t^2 u(t,x)=G(0)\partial_x^2 u(t,x)+\int_0^t G'(t-\tau)\partial_x^2 u(\tau,x)\,d\tau,    
\end{equation}
where we used the notation, $G'(t):=\partial_t G(t)$. Substituting the expression for the relaxation function, yields
\begin{equation}\label{motion2}
\partial_t^2 u(t,x)=D\,\partial_x^2 u(t,x)-\sum_{i=1}^N b_i\int_0^t e^{-r_i(t-\tau)}\partial_x^2 u(\tau,x)\,d\tau,    
\end{equation}
where $D$ is defined by
\begin{equation}\label{D}D=\sum_{i=1}^N s_i\,\,\text{with}\,\,b_i=s_i\,r_i,\,\,1\le i\le N.
\end{equation}
We omit specifying where the equations are satisfied unless it is unclear from the context. 

Equation \eqref{motion2} is nothing but the equation of motion for one of the well-known spring-dashpot models, called the extended Burgers model (abbreviated by EBM). To guarantee that the solutions of \eqref{IBP} decay exponentially as $t\rightarrow\infty$, we need to connect a spring in parallel to this model. In that case, we have $D>\sum_{i=1}^N b_i/r_i$ (see \cite{DKLNT, DLN}). Based on this, although it is a bit confusing in notations, we fix $b_i, r_i$'s and vary $D$. Hence, in this paper, we will consider all the following cases,
\begin{equation}  \label{cond_D}
D>h,\,\,D=h\ \text{and}\ D<h\,\,\text{with}\,\, h:=\sum_{i=1}^N b_i/r_i.
\end{equation} 
 
For simplicity, let $\xi u:=\partial_x^2 u $ and $I(r):=\int_0^t e^{-r(t-\tau)} \partial_x^2 u(\tau,x) \, d\tau $. Then, \eqref{motion2} takes the form,
\begin{equation}\label{I1.7}
\begin{aligned}
u''-D\xi u +  \sum_{i=1}^N b_iI(r_i)=0.
\end{aligned}
\end{equation}
% where $u'':=(u')'$ and “$'$" denotes the $t$-derivative.
Note that 
\begin{equation}\label{I1.8}
\begin{aligned}
I'(r)=\xi u-r I(r)
\end{aligned}
\end{equation}
holds. Then, introducing new dependent variables, $v=u'$ and $w_i=I(r_i)$ for $1\leq i\leq N$, we can write \eqref{I1.7} in the form of a system of partial differential equations given as
\begin{equation}\label{I1.9}
\left\{
\begin{array}{rcl}
u'&=&v,\\
v'&=&D\xi u +  \sum_{i=1}^N b_iw_i,\\
w_i'&=&\xi u -r_iw_i, \qquad\quad\ 1\leq i \leq N.
\end{array}
\right.
\end{equation}
We refer to this system as the augmented system. The matrix-vector representation of this system is given by
\begin{equation}\label{I1.10}
\begin{aligned}
U'=A_{N+2}U,
\end{aligned}
\end{equation}
where $A_{N+2}$ is a $N+2$ square matrix of the second-order partial differential operators in $x$ and $U=(u,v,w_1,\cdots,w_n)^{\mathfrak{t}}$ with $\mathfrak{t}$ denoting the transpose.

Consider $A_{N+2}$ as a densely defined closed operator on the $N+2$ number of products of $L^2(\Omega)$, denoted by $L^2(\Omega)^{N+2}$, with the domain
\begin{equation}\label{domain of A}
D(A_{N+2}):=\left\{
\begin{array}{ll}
U=(u,v,w_1,\cdots,w_N)^{\mathfrak{t}}\in H^2(\Omega)\times H^1(\Omega)\times L^2(\Omega)^N:\\
\qquad\qquad\qquad u(t,0)=\partial_x u(t,\pi/2)=0,\,t>0
\end{array}
\right\},
\end{equation}
where $H^s(\Omega),\,s=1,2$ are the $L^2(\Omega)$-based  Sobolev spaces of order $s=1,2$. Then, it is natural to analyze the spectrum of $A_{N+2}$, specifically, the eigenvalue problem:
\begin{equation}\label{I1.12}
\begin{aligned}
\lambda U=A_{N+2}\,U.
\end{aligned}
\end{equation}

A comprehensive numerical study on giving an approximate correspondence between the stretched exponential function and its approximation by the Prony series is given in \cite{MM}; here, it is noted that spectral information of $A_{N+2}$ with $D=h$ can be used for estimation of this function from data. General data analysis for the estimation of the stretched exponential function can be found in \cite{Phillips, Song} and the references therein, where the data can be spectral data obtained by the dynamical mechanical analysis (abbreviated by DMA) instruments.

\begin{comment}There is a lot of data analysis on the appropriate choice of $\beta,\,\tau$ of the stretched exponential function (see \cite{Phillips}, \cite{Song}, and the references therein).  However, it is very hard to find a paper on the eigenvalue problem \eqref{n1.7} other than \cite{LS}, which considered the case $D=h$. Needless to say that spectral information of $A_{N+2}$ with $D=h$ can be used for a data analysis of the stretched exponential function via the work on \cite{MM}. For instance, the spectral information can be obtained as measured data by using the DMA (dynamic mechanical analyzer).
\end{comment}

Due to the difficulty that comes from $A_{N+2}$ being non-self-adjoint, the authors of \cite{LS} considered clustered eigenvalues; we refer to this method as the clustered eigenvalues method. This method already appeared before in the inverse spectral geometry problem (see \cite{de Verdier}, \cite{Guillemin}, \cite{Gurarie}, \cite{Weinstein}). In our previous paper [1], we have given an inversion method for identifying $D$ and the relaxation tensor given as the Prony series by knowing two clustered sets of eigenvalues of $A_{N+2}$. As a natural continuation of our study, this paper aims to provide the numerical performance of our inversion. Combining this paper with our previous one [1] completes the spectral inversion method of the stress relaxation function in the glassy state, initiated from the articles \cite{LS} and \cite{MM}. We believe it to become a powerful data analysis method for this function.

As a consequence, the method reveals that it is a powerful method of data analysis, analyzing the relaxation of the glassy state data.  

\begin{comment}
Although we consider a bounded interval in our analysis, we note that it can be applied to the case $\Omega \subset {\mathbb R}^d$ with $d > 1$ being a bounded domain with a Lipschitz smooth boundary, and $\partial_x^2$ is replaced by a second-order positive operator defined in $L^2(\Omega)$. This is because the operator has an eigenfunction expansion and its eigenvalues $\{\lambda_\ell\}_{\ell=1}^\infty$ increasingly labeled, counting multiplicity and satisfying $\lambda_\ell \sim c_0\,\ell^{d/2}$ for $\ell \gg 1$, where $c_0>0$ is a constant independent of $\ell$ \cite{Agr, CH}. 
\end{comment}

The remainder of this paper is organized as follows. In Section 2, we recall the clustered eigenvalue problem more clearly and give a brief summary of our inversion method. Then, in Section 3, we show the achievement to our aim. Finally, we provide conclusions and discussions to our achievement.

\begin{comment}
Finally, we provide the numerical performance of our introduce the reduced eigenvalue problem with a parameter $k\in{\mathbb N}$ while describing the method mentioned above. The method can transform \eqref{n1.7} into an eigenvalue problem for a matrix $A_{N+2}^k$ which gives clustered eigenvalues of $A_{N+2}$. Then, in Section 3, we give the characterization of roots for the characteristic equation of $A_{N+2}^k$ and its limiting counterpart as $k\rightarrow\infty$. In Sections 4, 5, we present the speed of convergence of eigenvalues of $A_{N+2}^k$ as $k\rightarrow\infty$. In Section 6, we provide a numerical verification of the theoretical studies. In Section 7, we solve an inverse spectral problem to recover $D$ and the relaxation function from two clusters of eigenvalues. Finally, in the last section, we discuss our results and broader applications.
\end{comment}

\section{Inverse clustered eigenvalue problem and its inversion method}
%\noindent

In this section, we first explain the clustered eigenvalue problem for $A_{N+2}$, and its inverse problem. Then, we provide a summary of our inversion method for this inverse problem. To begin studying the eigenvalue problem for $A_{N+2}$, we consider an eigenvalue problem as follows. Namely, we replace $\xi$ with its generic eigenvalue $-(2k-1)^2$ with $k\in{\mathbb N}$. More precisely, we search for an eigenvalue $\eta$ of $-\partial_x^2$ which gives a non-trivial solution to the following boundary value problem:
\begin{equation}\label{n2.2}
\left\{
\begin{array}{ll}
\partial_x^2u-\eta u=0,\,\, x\in (0,\frac{\pi}{2}),\\
u(0)=0,\,\, \partial_xu(\frac{\pi}{2})=0.
\end{array}
\right.
\end{equation}
It is easy to see that a possible eigenvalue $\eta$ and its associated eigenvector $u$ are given as
\begin{equation}\label{n2.3}
\begin{aligned}
\eta=-(2k-1)^2,\, u=\sin (2k-1)x, \quad k\in \mathbf{N}.
\end{aligned}
\end{equation}

Now, let us see that this eigenvalue problem can generate a cluster of eigenvalues of $A_{N+2}$ as follows. Denote
$A_{N+2}^k$ as $A_{N+2}$ with $\xi$ replaced by $-(2k-1)^2$. Then $A_{N+2}^k$ is no longer a differential operator. It is just a multiplication operator on ${\mathbb R}^{N+2}$ by this $(N+2)\times (N+2)$ real matrix $A_{N+2}^k$.
Then, consider the following {\sl clustered eigenvalue problem}:
\begin{equation}\label{new1.7}
\begin{aligned}
\lambda\,\tilde{U}^k=A_{N+2}^k\,\tilde{U}^k,
\end{aligned}
\end{equation}
where $\tilde{U}^k\in{\mathbb C}^{N+2}$. This eigenvalue $\lambda$ with $\lambda\not=-r_i,\,1\le i\le N$ is an eigenvalue of \eqref{I1.12}, and its associated eigenfunction $U^k$ is given as  $U^k=(u,v,w_1,\cdots,w_n)^{\mathfrak{t}}$ of \eqref{I1.12} given as 
\begin{equation}\label{n2.4}
\left\{
\begin{array}{ll}
u=\sin (2k-1)x,\\
v=\lambda\sin (2k-1)x,\\
w_i=-(\lambda+r_i)^{-1}(2k-1)^2\sin (2k-1)x, \quad 1\leq i\leq N.
\end{array}
\right.
\end{equation}
For each $k\in{\mathbb N}$, the eigenvalues of \eqref{new1.7} are called the {\sl clustered eigenvalues}. 

\medskip
Then, the inverse problem associated with the clustered eigenvalue problem is stated as follows.

\medskip
\noindent
{\bf Inverse problem:} Recover $D,\, b_i, r_i$ with $1\le i\le N$ of \eqref{motion2} by knowing sets of the clustered eigenvalues of \eqref{new1.7}. Here, we emphasize that we consider all the cases \eqref{cond_D}.

\medskip
An answer to the inverse problem is stated as follows.
\begin{theorem}\label{thm_IP}
By knowing two clusters of eigenvalues associated with $k=k_1,\,k_2\in{\mathbb N}$, we can recover $D,\, b_i, r_i$ with $1\le i\le N$ of \eqref{motion2}.
\end{theorem}

Before giving a proof of this theorem, we give some preliminaries. We first note that from \eqref{I1.12}, we have 
\begin{equation}\label{n1.8}
\begin{aligned}
\lambda^2u+D(2k-1)^2 u +  \sum_{i=1}^N b_iw_i=0,
\end{aligned}
\end{equation}
and
\begin{equation}\label{n1.9}
\begin{aligned}
\lambda w_i=-(2k-1)^2 u -r_iw_i,\,\,1\le i\le N,
\end{aligned}
\end{equation}
which are the same to
\begin{equation}\label{n1.10}
\begin{aligned}
(\lambda+r_i)w_i=-(2k-1)^2 u,\,\,\,1\le i\le N.
\end{aligned}
\end{equation}
Multiplying $\Pi_{1\leq j\leq N}(\lambda+r_j)$ to \eqref{n1.8} and using \eqref{n1.10}, we obtain
\begin{equation}\label{n3.2}
\begin{aligned}
P_N^k(\lambda):=(D+\frac{\lambda^2}{(2k-1)^2})\Pi_{1\leq j\leq N}(\lambda+r_j) -  \sum_{i=1}^N b_i\Pi_{1\leq j\leq N,j\neq i}(\lambda+r_j) =0.
\end{aligned}
\end{equation}
We refer to \eqref{n3.2} and $P_N^k(\lambda)$ as the characteristic equation and the characteristic polynomial associated with $\partial_t-A_{N+2}^k$, respectively. From the arguments deriving \eqref{n3.2}, it is clear that for each $k\in{\mathbb N}$, the roots of \eqref{n3.2} are the eigenvalues of $A_{N+2}$. These are the preliminaries, and we can move to a proof.

\medskip
\noindent
{\bf Proof of Theorem \ref{thm_IP}:}
The idea of a proof is based on the well-known fact that for a monic polynomial, knowing all of its roots is equivalent to knowing the polynomial.
For further details of the proof, we need to quote Lemma 3.2 of \cite{CDDLN} given as follows:
\begin{lemma}\label{important}
There exists at least $N$ real roots of the characteristic polynomial $P_N^k(\lambda)$ such that 
\begin{equation}\label{n4.7}
\begin{aligned}
-r_N<a^k_N<-r_{N-1}<a^k_{N-1}<\cdots<-r_2<a^k_2<-r_1<a^k_1
\end{aligned}
\end{equation}
and
\begin{equation}
\label{expression_D2}
\sum_{i=1}^N \frac{b_i}{a^k_j+r_i} =D+\frac{(a_j^k)^2}{(2k-1)^2},\,\,1\le j\le N.
\end{equation}
The other two roots are contained in the set $B^D_-\cup \{(-r_N,a_1^k)\} $, where $B^D_-:=\{c+id :\frac{-r_N-a_1^k}{2}<c<\frac{-r_1-a_1^k}{2}\} $.
\end{lemma}

Based on this lemma, for $k_1, k_2\in{\mathbb N},\,k_1<k_2$, assume that we know the clusters of eigenvalues associated to $k_1,\,k_2$ given as 
\begin{equation}\label{i7.1}
\begin{aligned}
a^{k_1}_1, a^{k_1}_2, a^{k_1}_3, \cdots,a^{k_1}_{N}, a^{k_1}_{N+1}=p^{k_1}+iq^{k_1},  a^{k_1}_{N+2}=p^{k_1}-iq^{k_1}
\end{aligned}
\end{equation}
and
\begin{equation}\label{i7.2}
\begin{aligned}
a^{k_2}_1, a^{k_2}_2, a^{k_2}_3, \cdots,a^{k_2}_{N}, a^{k_2}_{N+1}=p^{k_2}+iq^{k_2},  a^{k_2}_{N+2}=p^{k_2}-iq^{k_2},
\end{aligned}
\end{equation}
where $a^{k_1}_{N+1}$, $a^{k_1}_{N+2}$ or $a^{k_2}_{N+1}$, $a^{k_2}_{N+2}$ may be real.

\medskip
Let us recall the characteristic polynomial $P^k_N(\lambda)$ given by \eqref{n3.2}. Then, since
this polynomial is of degree $N+2$ and its coefficient of $\lambda^{N+2}$ is $\frac{1}{(2k-1)^2}$, we have 
\begin{equation}\label{i7.4}
\begin{aligned}
P^{k_1}_N(\lambda)=\frac{1}{(2k_1-1)^2}\cdot\Pi_{1\leq j\leq N+2}(\lambda-a_j^{k_1})
\end{aligned}
\end{equation}
and
\begin{equation}\label{i7.5}
\begin{aligned}
P^{k_2}_N(\lambda)=\frac{1}{(2k_2-1)^2}\cdot\Pi_{1\leq j\leq N+2}(\lambda-a_j^{k_2})
\end{aligned}
\end{equation}
from \eqref{i7.1} and \eqref{i7.2}, respectively.
By putting $k=k_1$ and $k=k_2$ in \eqref{n3.2}, we have
\begin{equation}\label{i7.6}
\begin{aligned}
P^{k_1}_N(\lambda)-P^{k_2}_N(\lambda)&=\frac{1}{(2k_1-1)^2}\lambda^2\Pi_{1\leq j\leq N}(\lambda+r_j)-\frac{1}{(2k_2-1)^2}\lambda^2\Pi_{1\leq j\leq N}(\lambda+r_j)\\
&=[\frac{1}{(2k_1-1)^2}-\frac{1}{(2k_2-1)^2}]\cdot\lambda^2\Pi_{1\leq j\leq N}(\lambda+r_j).
\end{aligned}
\end{equation}
This implies
\begin{equation}\label{i7.7}
\begin{aligned}
\Pi_{1\leq j\leq N}(\lambda+r_j)=[\frac{1}{(2k_1-1)^2}-\frac{1}{(2k_2-1)^2}]^{-1}\lambda^{-2}\big(P^{k_1}_N(\lambda)-P^{k_2}_N(\lambda)\big).
\end{aligned}
\end{equation}
Here, note that we know $k_1, k_2$ and all the respective roots of \eqref{i7.4} and \eqref{i7.5}, the polynomial on the right-hand side
of \eqref{i7.7} is known due to the mentioned well-known fact. Hence we know $r_j,\,1\leq j\leq N$.

Observe that for each $1\leq i \leq N$, we have
\begin{equation}\label{n4.8}
\begin{aligned}
P^k_N(-r_i)&= -   b_i\Pi_{1\leq j\leq N,j\neq i}(-r_i+r_j) \\
&=(-1)^i |b_i\Pi_{1\leq j\leq N,j\neq i}(-r_i+r_j)|.
\end{aligned}
\end{equation}

Using this for $k=k_2$, we can recover $b_i$
for $1\leq i\leq N$ because we know $P_N^{k_2}(\lambda)$ and $r_i,\,1\leq i\leq N$. We can further recover $D$ using 
\eqref{n3.2} with $k=k_2$ because we already know $b_i,\,r_i,\,1\le i\le N$. Thus, we have recovered $D$, $r_i$ and $b_i$ for $1\leq i\leq N$ by knowing the eigenvalues of $P^{k_1}_N(\lambda)$ and $P^{k_2}_N(\lambda)$ with $k_1, k_2 \in \mathbb{N}$, $k_1<k_2$. 
Thus, we have completed the proof.

\section{Numerical performance of the inversion method}
%\noindent

In this section, we present several numerical results to demonstrate the numerical performance of our inversion method, given in Theorem \ref{thm_IP}. We will use a bisection method to realize the mentioned method. For applying the bisection method, we point out the following important considerations.

\medskip\noindent
{\bf Considerations for numerically verifying the inversion method:} 

\medskip\noindent
Let
\begin{equation}\label{eq_Q}
Q(\lambda;k_1,k_2):=[\frac{1}{(2k_1-1)^2}-\frac{1}{(2k_2-1)^2}]^{-1}\lambda^{-2}\big(P^{k_1}_N(\lambda)-P^{k_2}_N(\lambda)\big).
\end{equation} 
For each fixed $\ell=1,2$, $Q(\lambda;k_1,k_2)$ has a single root $-r_{j-1}$ in each $(a_j^{k_\ell},a_{j-1}^{k_\ell})$ with $j=2,\cdots,N$ by the Lemma in Section 2 and \eqref{i7.7}. 
This implies that the condition $Q(a_j^{k_\ell};k_1,k_2)\,Q(a_{j-1}^{k_\ell};k_1,k_2)<0$ to start the bisection method holds, enabling us
to get $-r_{j-1}$. Upon having $-r_{j-1},\,j=2,\cdots, N$, the remaining $-r_N$ can be easily obtained. Further, by the argument of proof of Theorem \ref{thm_IP}, $b_j$'s and $D$ can be easily obtained using \eqref{n4.8} and \eqref{n3.2}, respectively. 

\medskip
Now, it is worth mentioning that in order to shorten the DMA measurement time and improve the approximation accuracy near the origin of the stretched exponential function approximated by the Prony series frequency, we set $k=k_1, k_2$ fairly large within the DMA frequency range. More precisely, we take $k\in[81,1001]$. Referring to \cite{MM}, the number $N$ of terms in the Prony series is set $N=5$ or $N=9$. Then, giving the target parameters $D,\,b_i,\,r_i$ with $1\le i\le N$, we generate the measured data, two sets of clustered eigenvalues $\{a_j^k\}_{j=1}^{N+2}$ with $k=k_1, k_2$, following the numerical study \cite{CDDLN} of the clustered eigenvalue problem. Here, by setting $h:=\sum_{i=1}^N {b_i}/{r_i}=1$, we take one of the target parameter $D$ as $D=0.5, 1, 5$ to satisfy all the cases $D>h,\,D=h,\,D<h$. Further, we take the other target parameters $\{b_i\}_{i=1}^N$ and $\{r_i\}_{i=1}^N$ as $b_i=1,\,r_i=5\,i,\,i=1,2,\cdots N$.

\begin{comment}
In order to get the eigenvalues as measurements of the inverse problem, we compute the roots, $a_1^k,\,a_2^k,\,\cdots,\,a_N^k$, of \eqref{n3.2} for two different $k$. More precisely, by fixing $h,\,N,\,\text{$r_i$'s}$ to $h=1$, $N=5\,\,\text{or}\,\,9$, $r_i=5*i,\,i=1,\,2,\,\cdots,\,N$, we numerically calculate the roots $a_j^k$ of $P_N^k(\lambda)$ defined by \eqref{n3.2} for different $k, N$ for the cases $D>h,\,D=h,\,D<h$, i.e., $D=0.5,\,1,\,5$. 

From the proof of Theorem \ref{thm_IP}, our numerical reconstruction can be decomposed into three steps. First, we apply the bisection method to recover $r_i, 1\leq i\leq N$ by using \eqref{i7.7}. In fact, the parameters $-r_i,\,i=1,\,2,\,\cdots,\,N$ make the left side of \eqref{i7.7} vanish. Recovering $r_i$'s can be transformed into searching the roots of the right side of \eqref{i7.7}, which are $-r_i,\,i=1,\,2,\,\cdots,\,N$ satisfying \eqref{n4.7}. Second, formula \eqref{n4.8} can be used directly to recover $b_i$, $i=1,2, \ldots, N$. Finally, one can use \eqref{n3.2} to recover $D$.
\end{comment}

For all the numerical results given below after References, the noise in the measurements is given by perturbing $a_N^k$ as follows
\begin{equation}\label{noisydata_add}
a_N^{k,\delta}:=\left(1+\delta\times(2\times {\rm rand(1)} -1)\right)\times a_N^k,
\end{equation}
where $\delta$ is the relative noise level, and $\rm rand(1)$ generates a uniformly distributed random number in the interval $(0,\,1)$.
 
\begin{figure}[htp]
\centering
\begin{tabular}{ccc}
 {\small $(k_1,\,k_2)=(81,\,91)$} & {\small $(k_1,\,k_2)=(81,\,501)$} & {\small $(k_1,\,k_2)=(81,\,1001)$ }\\
\adjustbox{trim=6.3mm 0mm 5.3mm 2.5mm, clip}{
\includegraphics[width=0.4\textwidth]{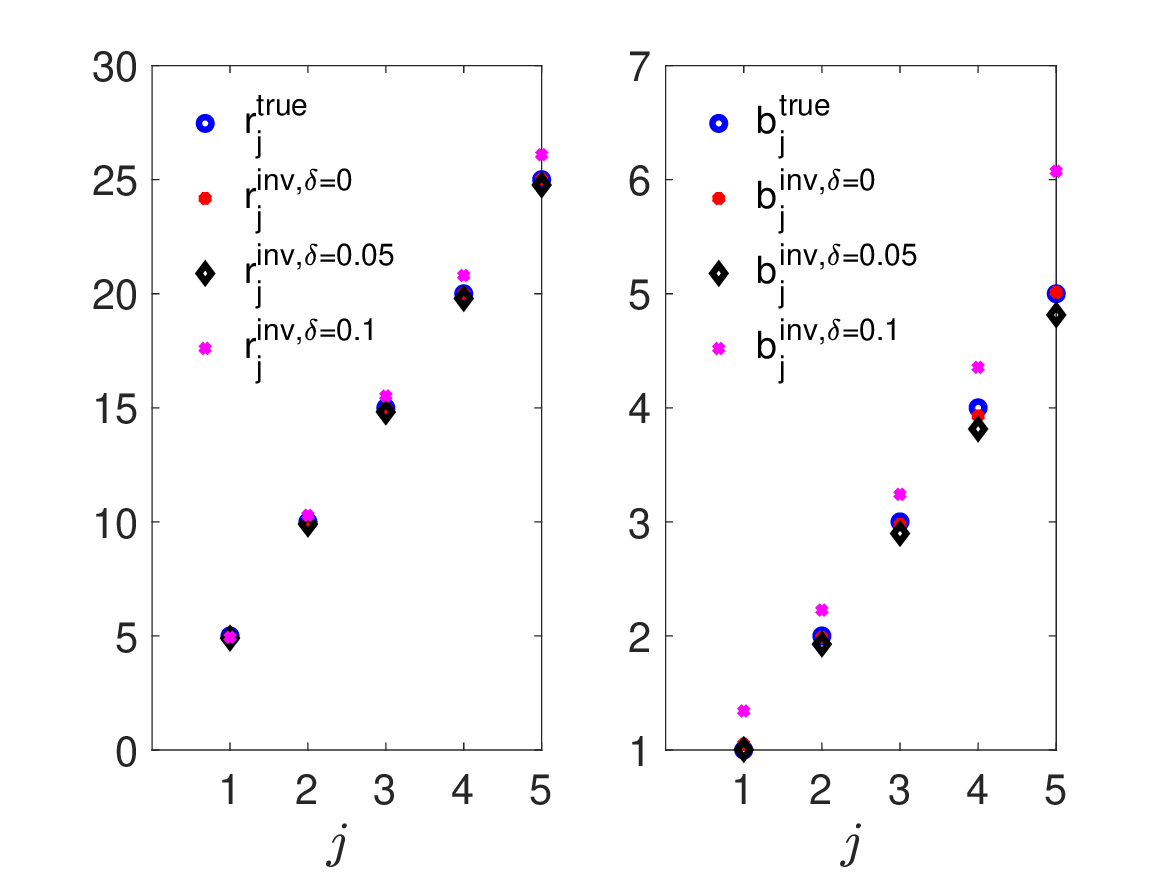}}
&\adjustbox{trim=4.9mm 0mm 5.3mm 2.5mm, clip}{\includegraphics[width=0.4\textwidth]{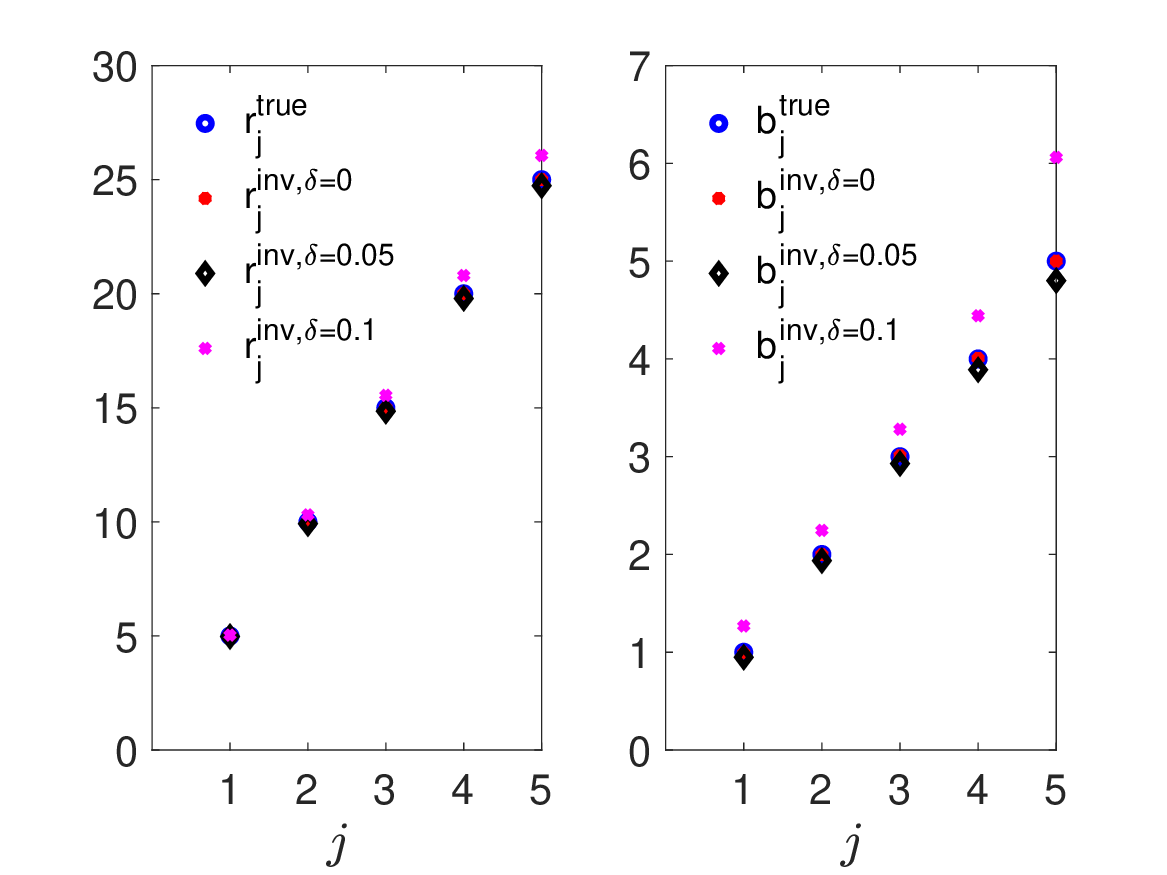}}
&\adjustbox{trim=4.9mm 0mm 5.3mm 2.5mm, clip}{\includegraphics[width=0.4\textwidth]{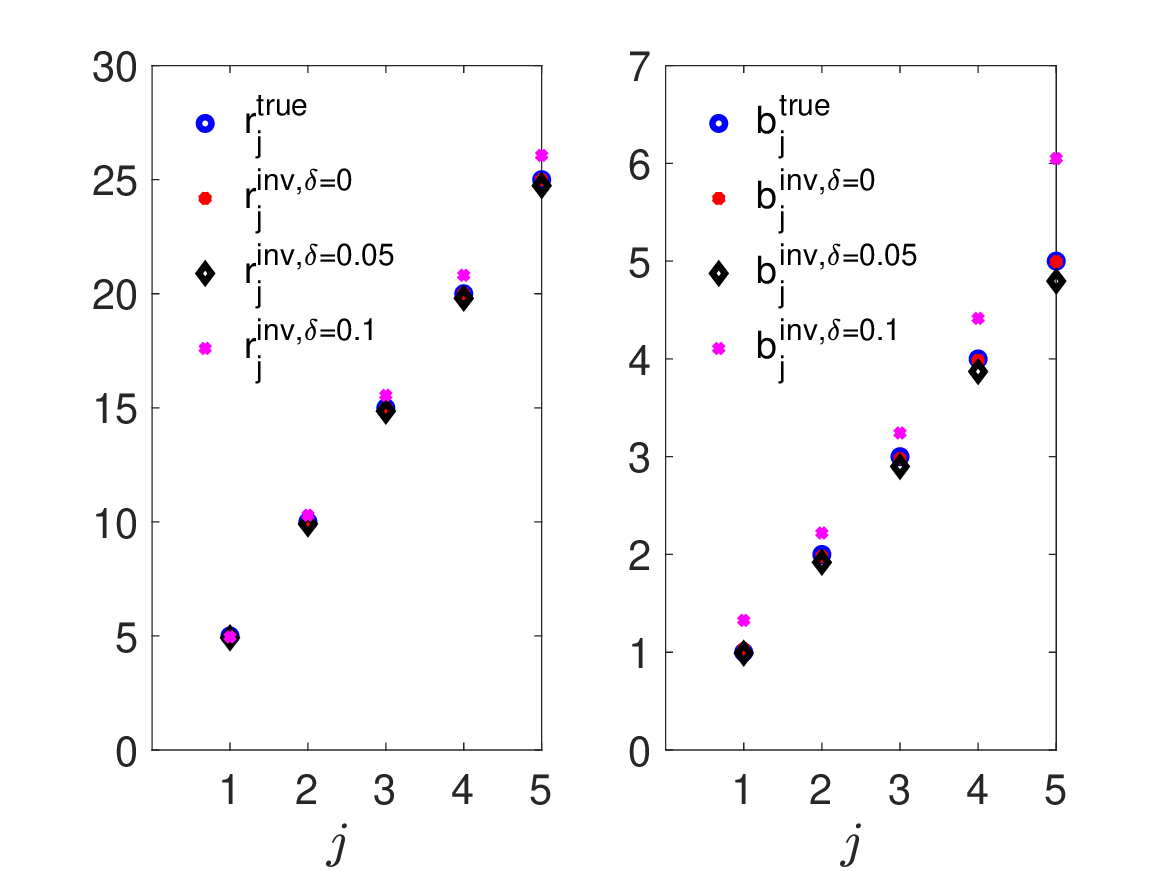}}\\
\adjustbox{trim=6.3mm 0mm 5.3mm 2.5mm, clip}{
\includegraphics[width=0.4\textwidth]{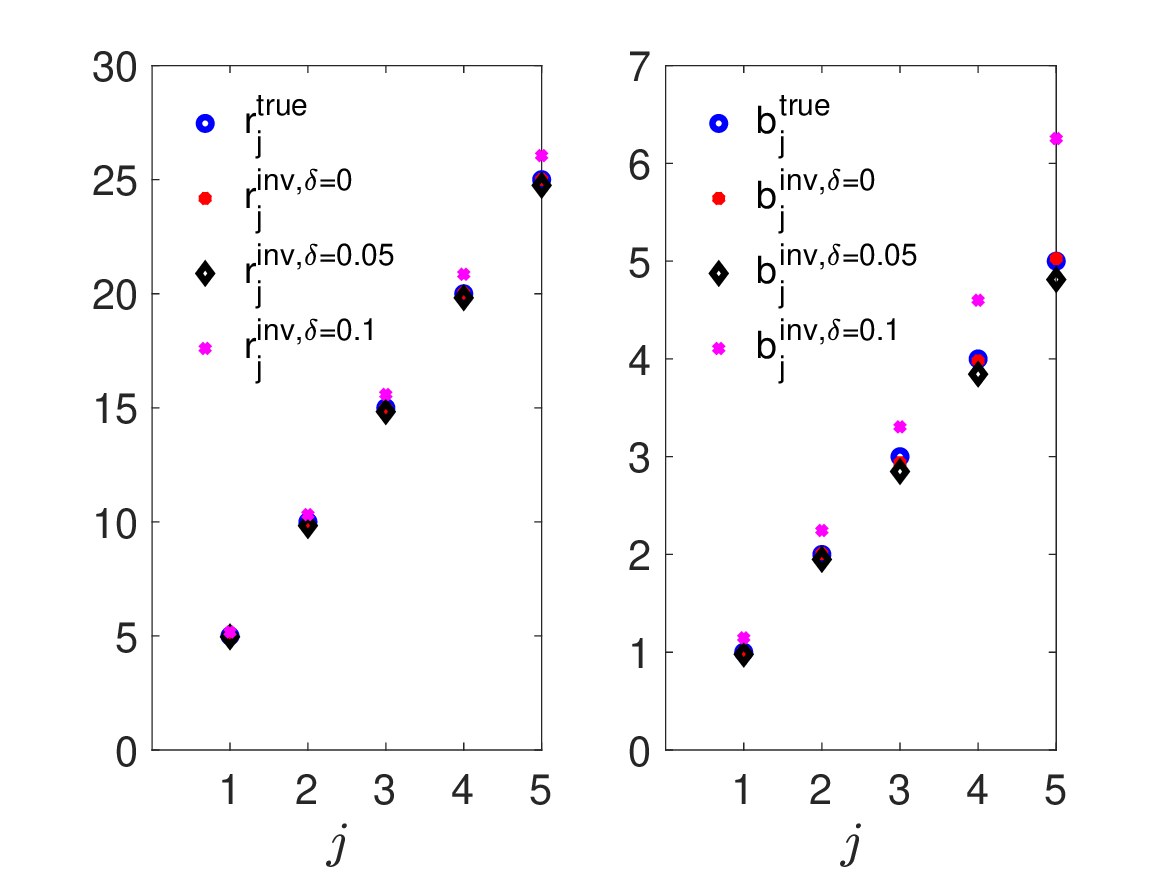}}
&\adjustbox{trim=6.3mm 0mm 5.3mm 2.5mm, clip}{
\includegraphics[width=0.4\textwidth]{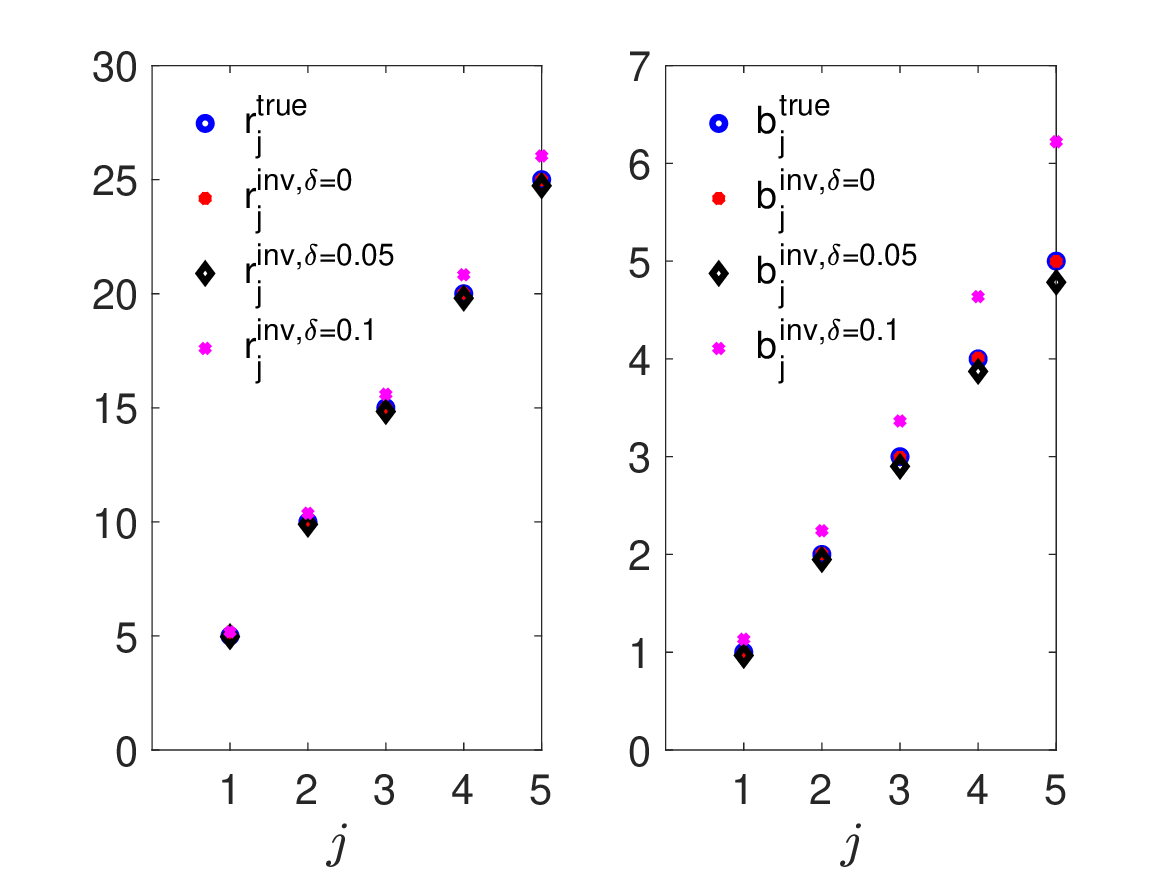}}
&\adjustbox{trim=6.3mm 0mm 5.3mm 2.5mm, clip}{
\includegraphics[width=0.4\textwidth]{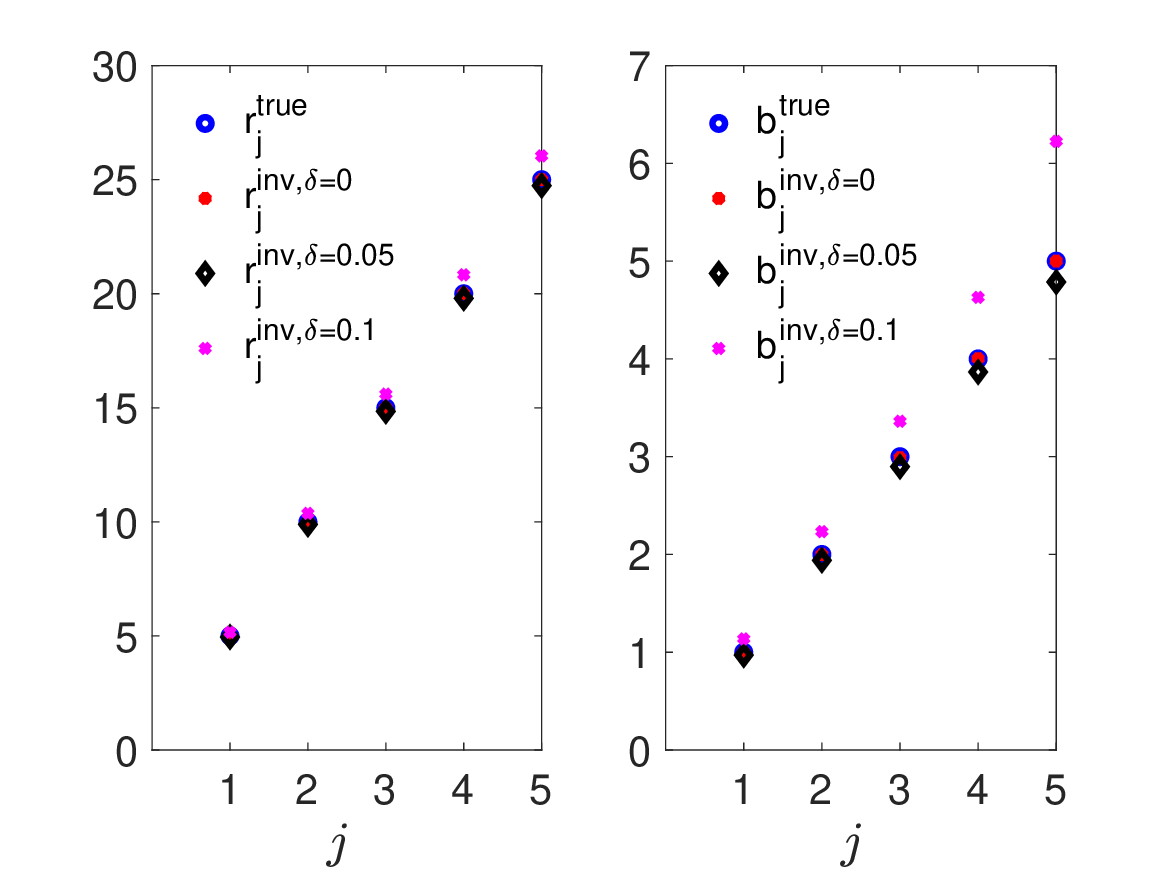}}\\
\adjustbox{trim=4.9mm 0mm 5.3mm 2.5mm, clip}{\includegraphics[width=0.4\textwidth]{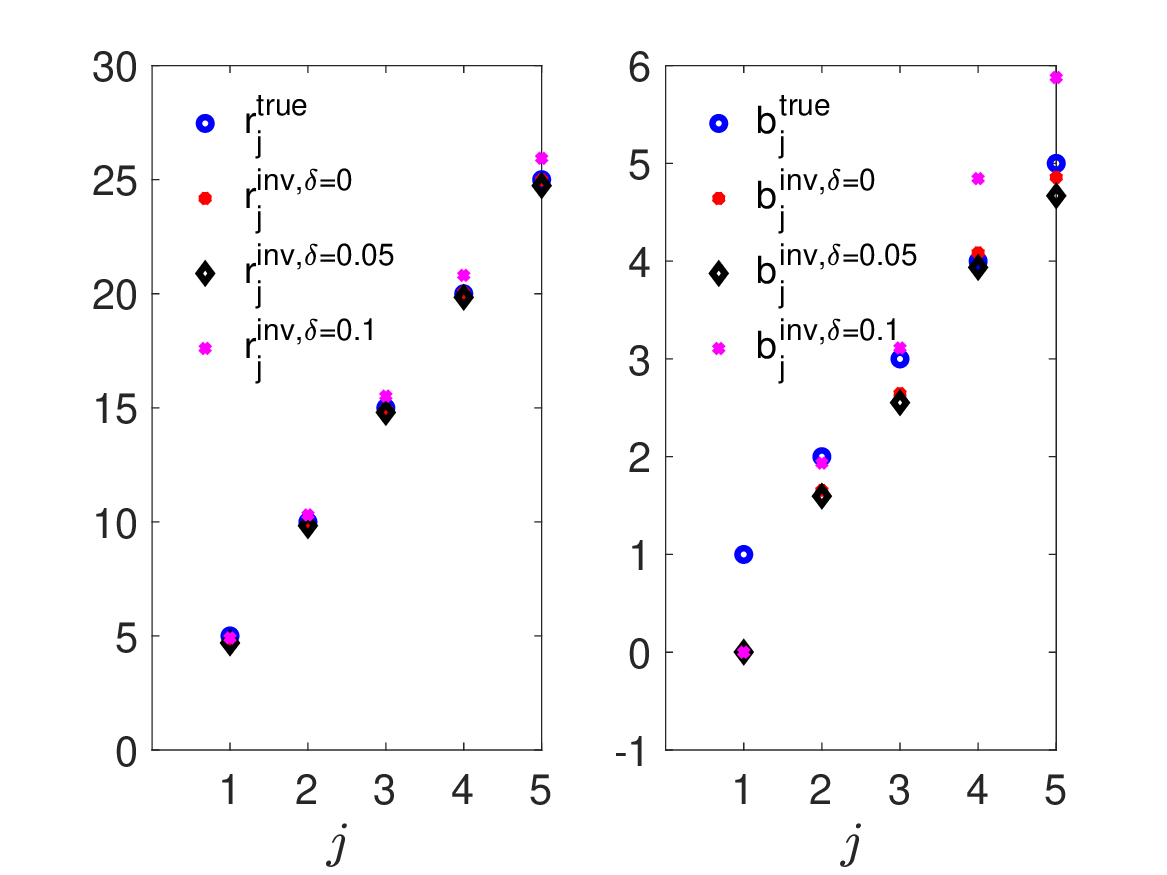}}
&\adjustbox{trim=4.9mm 0mm 5.3mm 2.5mm, clip}{\includegraphics[width=0.4\textwidth]{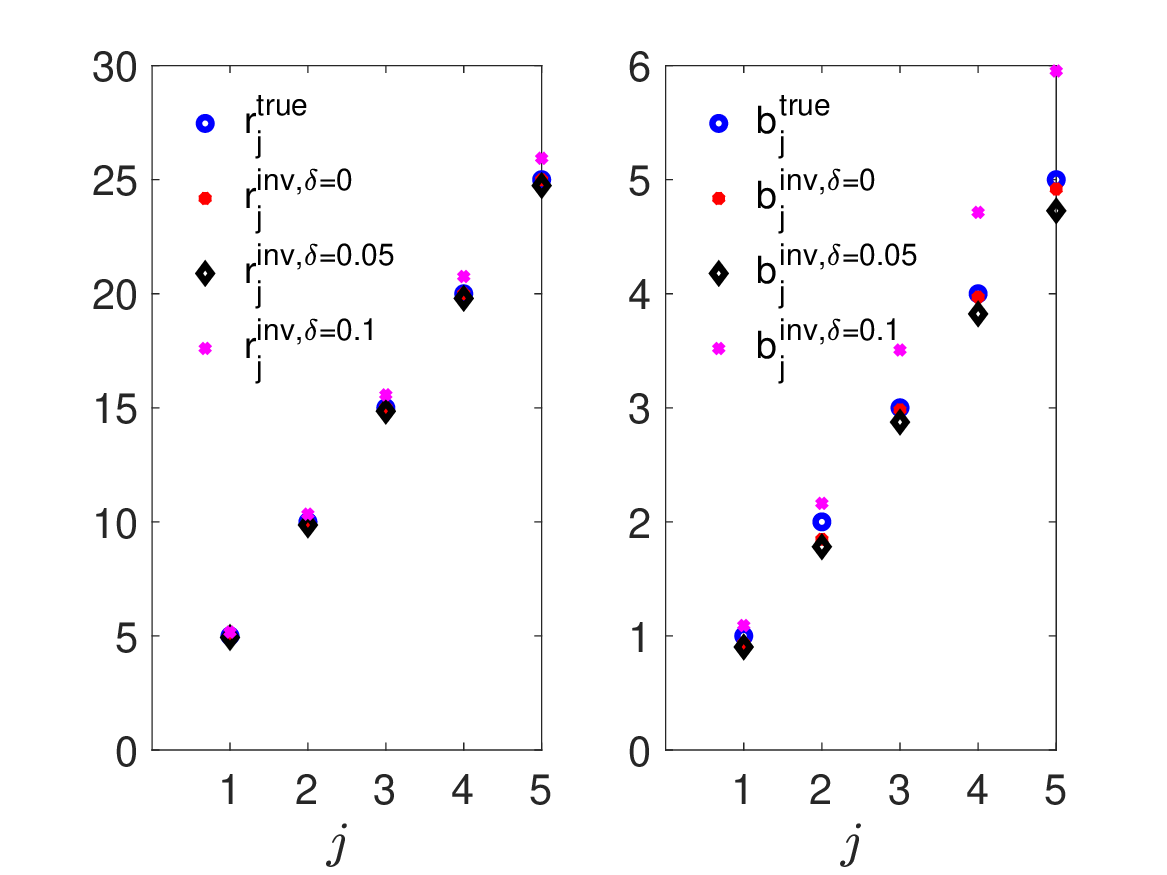}}
&\adjustbox{trim=4.9mm 0mm 5.3mm 2.5mm, clip}{\includegraphics[width=0.4\textwidth]{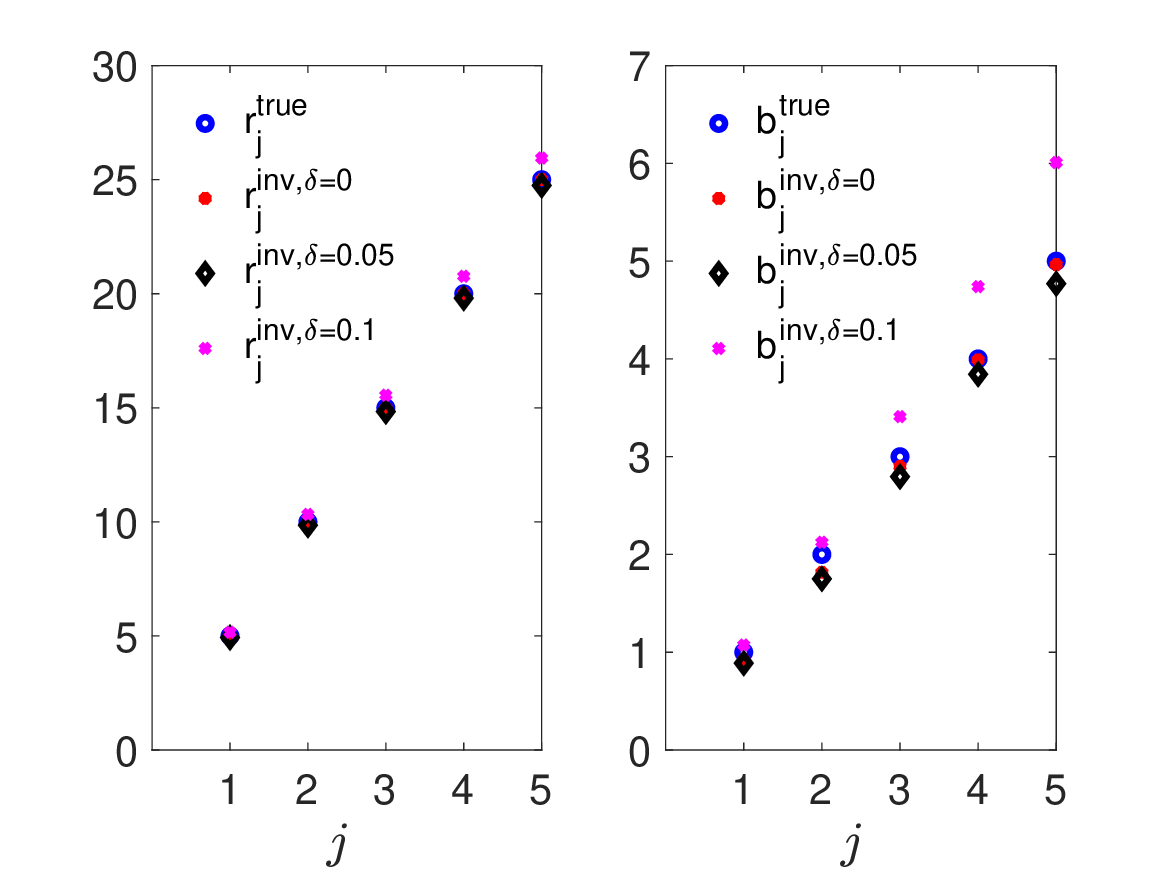}}
\end{tabular}
\caption{Comparisons between reconstructions $(r_j^{inv},\,b_j^{inv},\,D^{inv})$ and true values $(r_j^{true},\,b_j^{true},\,D^{true})$ for $N=5$ and noise level $\delta=0,\,0.05,\,0.1$. Here, each row represents the reconstructions for $D=0.5,\,1,\,5$, respectively. Each column represents the reconstructions for the pair of frequencies $(k_1,\,k_2)=(81,\,91),\,(81,\,501),\,(81,\,1001)$
}

\label{fig_N5}
\end{figure}

\begin{example}\label{exm_1}
Let $N=5$ and $D=0.5,\,1,\,5$. The reconstructions of $(r_j^{inv},\,b_j^{inv}),\,j=1,\,2,\,\cdots,\,N$ and $D^{inv}$ are plotted in Figure \ref{fig_N5} and listed in Table \ref{tbl:N5}, respectively. 
\end{example}

From the results shown in Figure \ref{fig_N5}, the reconstructions of $(r_j^{inv},b_j^{inv}),\,j=1,\,\cdot,\,N$ are indistinguishable compared with the exact solutions for the noise-free case. Also, from Table \ref{tbl:N5}, it is easy to observe that the reconstruction of $D^{inv}$ converges to the exact solution as the distance between $k_1$ and $k_2$ increases. That is, two clustered eigenvalues corresponding to $k_1$ and $k_2$ should be as different as possible. Based on this and from the convergence of the eigenvalues $a_j^k$ on $k$ shown in \cite{CDDLN}, we have to choose a lower frequency $k_1$ if $k_1< k_2$. To conserve space, we have omitted the presentation of the numerical results computed by using lower frequency eigenvalues.

For the noisy case, the reconstructions are affected by noise and become worse as the noise level increases. But they are still acceptable even for 10\% noise level. Comparing the reconstructions $r_j^{inv}$ and $b_j^{inv}$ with the exact solutions for each pair $(k_1,\,k_2)$, $r_j^{inv}$ is always better than $b_j^{inv}$.  This is thought to be because our inversion method is a step-by-step process in which $r_j$'s are reconstructed first, and then $b_j$'s are reconstructed, and hence reconstruction of $b_j$'s is affected to measurement errors than the reconstruction of $r_j$'s. Regarding the reconstruction of $D$ and comparing it with the noise-free case, we note that the reconstruction of $D_j^{inv}$ no longer converges to the true one as the distance between $k_1$ and $k_2$ increases.

\begin{figure}[htp]
\centering
\begin{tabular}{ccc}
 {\small $(k_1,\,k_2)=(81,\,91)$}  & {\small $(k_1,\,k_2)=(81,\,501)$} & {\small $(k_1,\,k_2)=(81,\,1001)$}\\
\adjustbox{trim=5.1mm 0mm 5.4mm 2.5mm, clip}{
\includegraphics[width=0.38\textwidth]{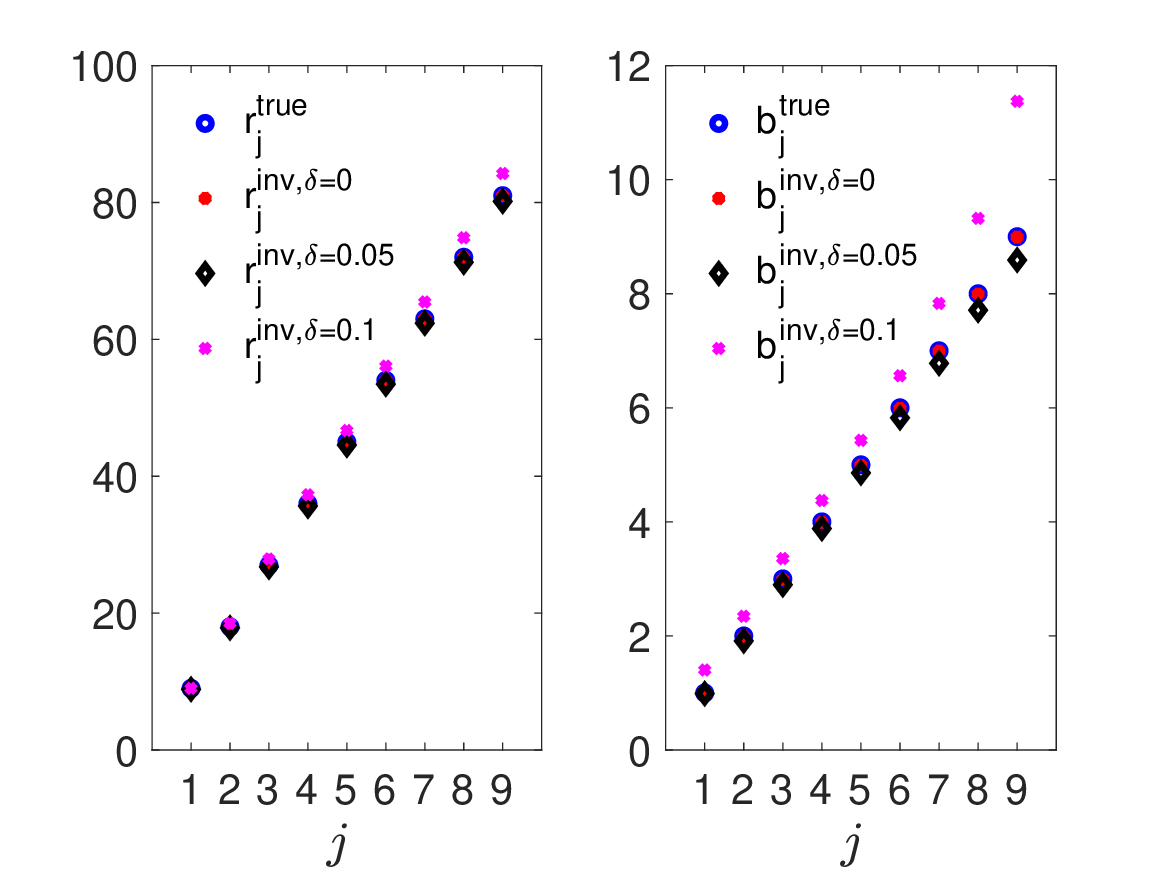}}
&\adjustbox{trim=3.7mm 0mm 5.4mm 2.5mm, clip}{\includegraphics[width=0.38\textwidth]{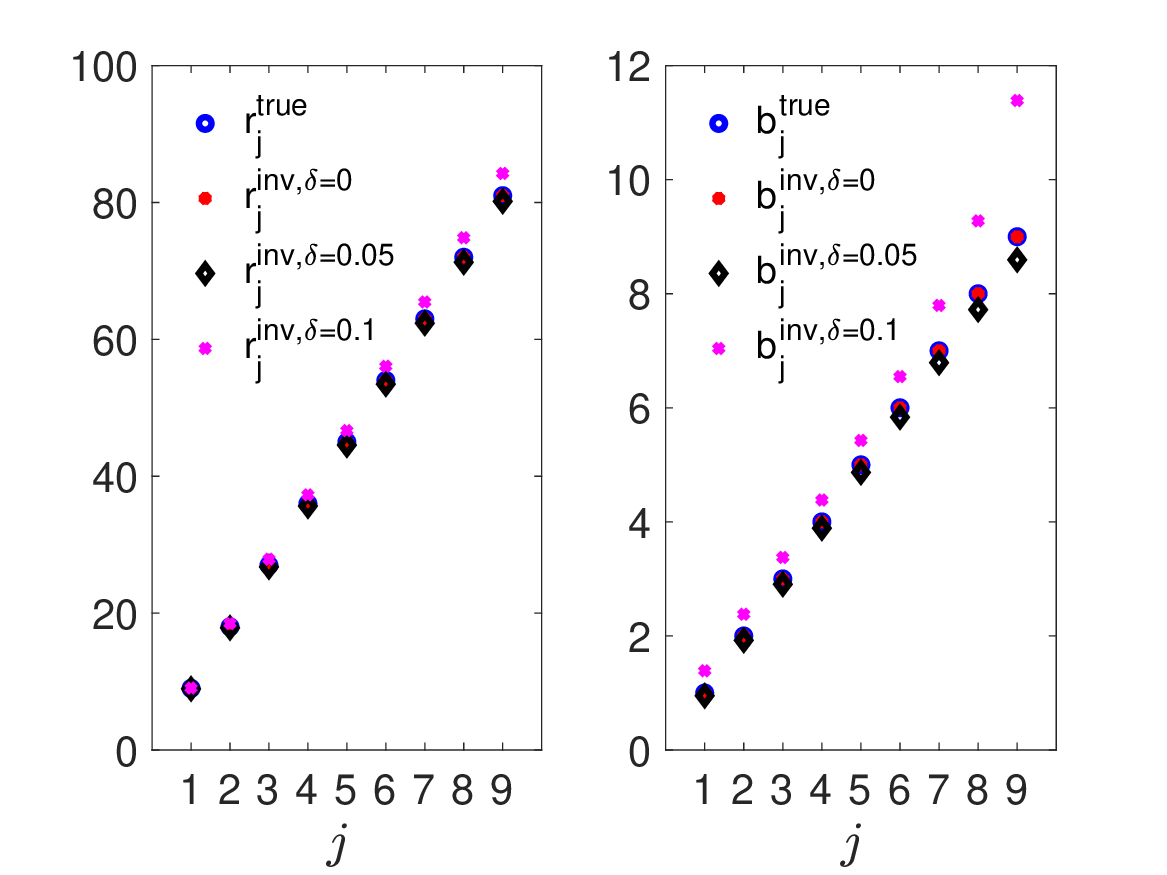}}
&\adjustbox{trim=3.7mm 0mm 5.4mm 2.5mm, clip}{\includegraphics[width=0.38\textwidth]{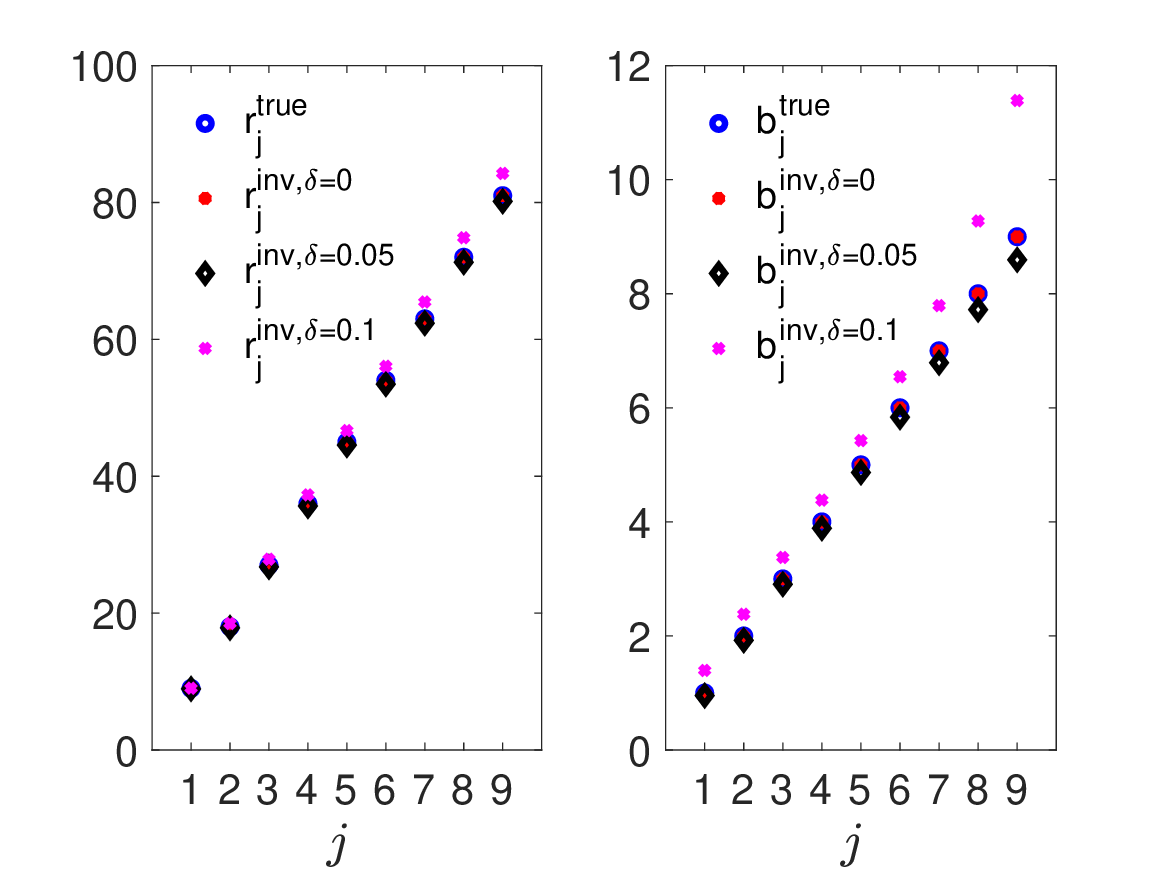}}\\
\adjustbox{trim=5.15mm 0mm 5.45mm 2.5mm, clip}{
\includegraphics[width=0.38\textwidth]{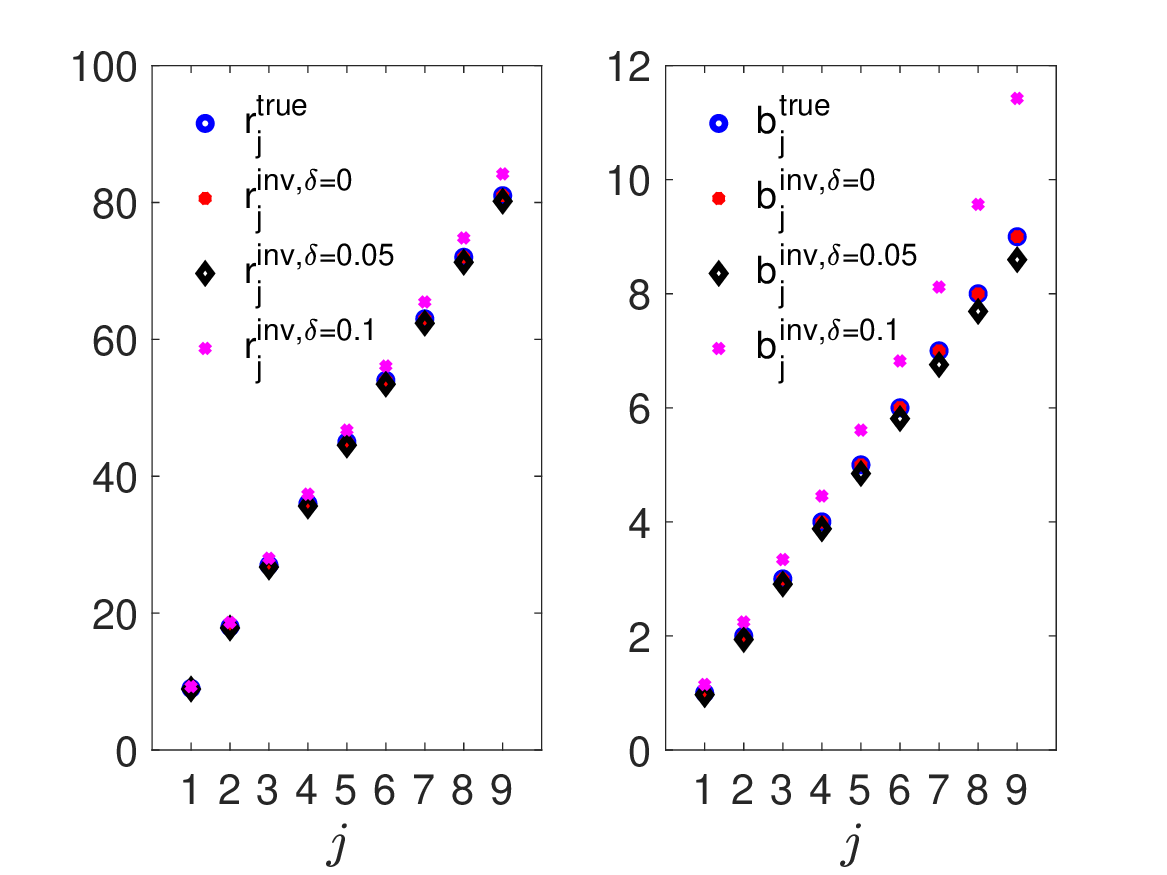}}
&\adjustbox{trim=3.75mm 0mm 5.4mm 2.5mm, clip}{
\includegraphics[width=0.38\textwidth]{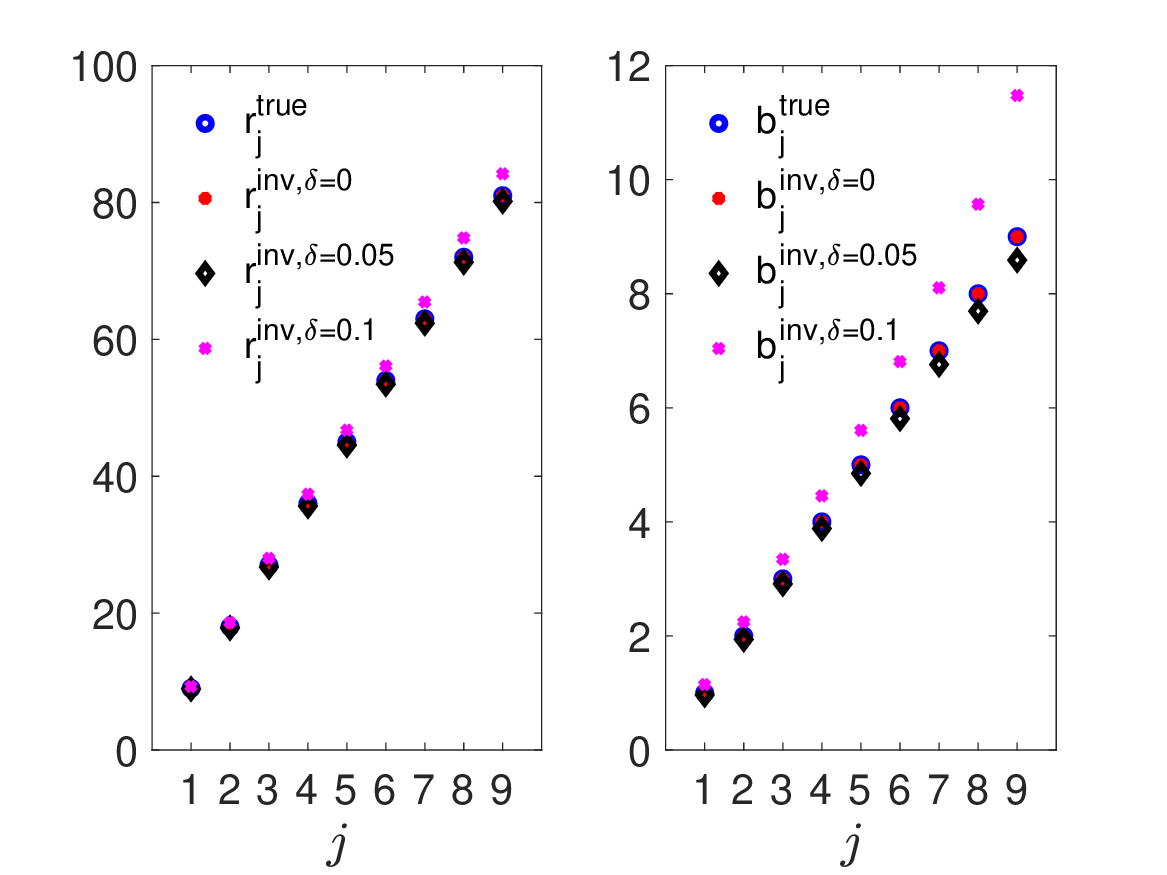}}
&\adjustbox{trim=3.7mm 0mm 5.4mm 2.5mm, clip}{
\includegraphics[width=0.38\textwidth]{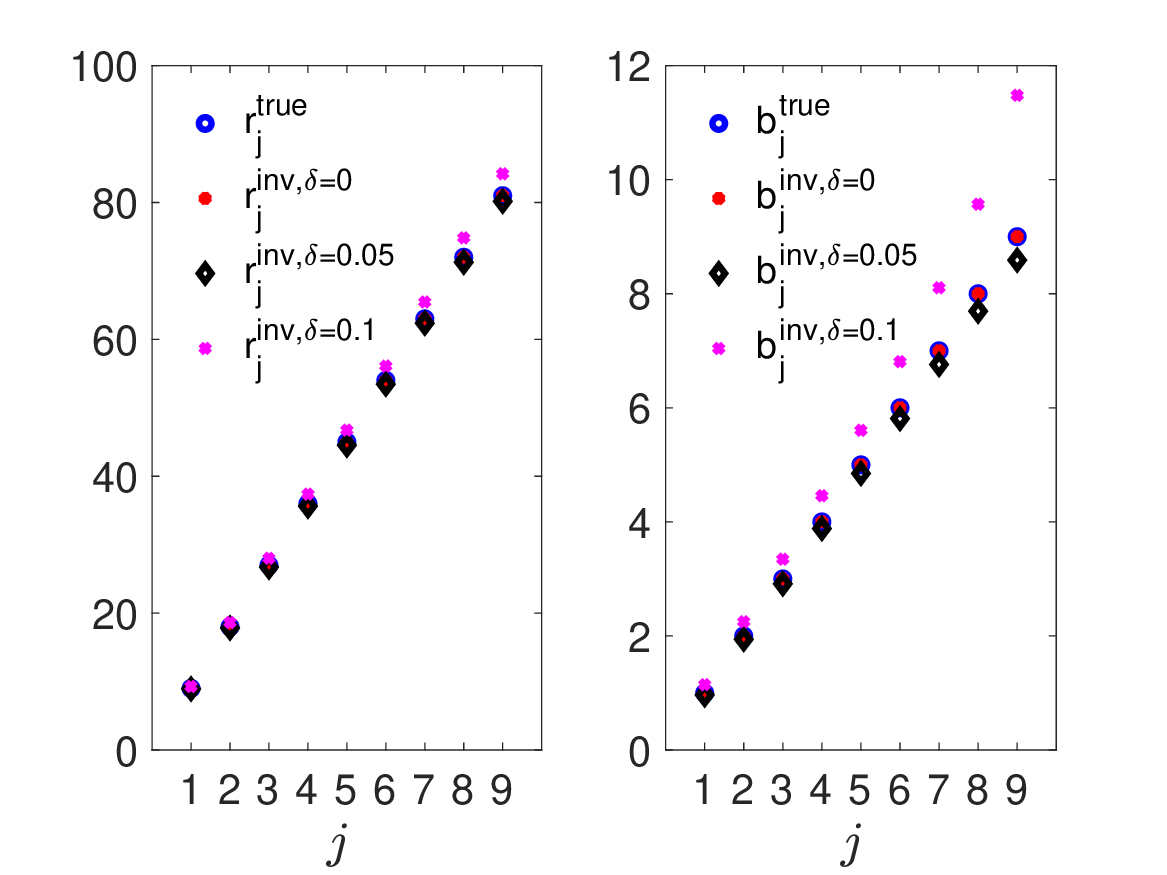}}\\
\adjustbox{trim=3.7mm 0mm 5.4mm 2.5mm, clip}{\includegraphics[width=0.38\textwidth]{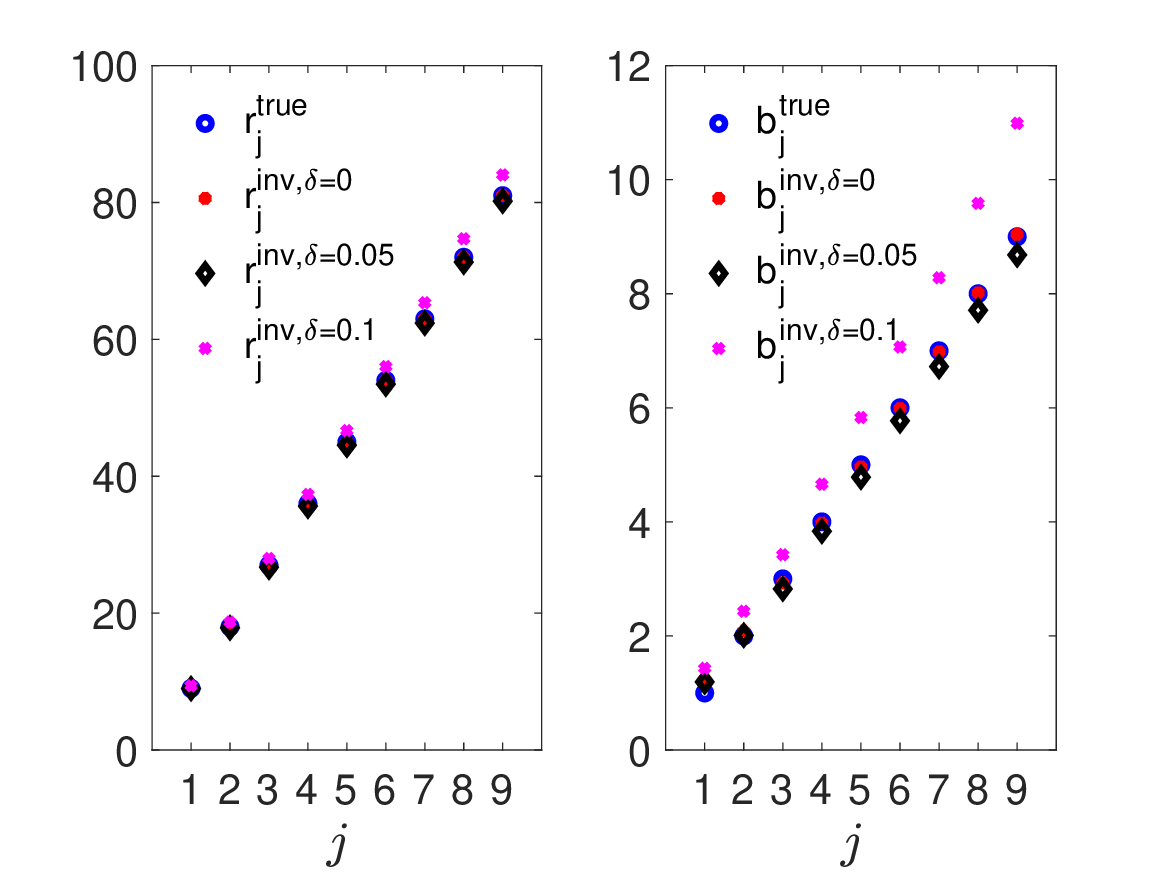}}
&\adjustbox{trim=3.7mm 0mm 5.4mm 2.5mm, clip}{\includegraphics[width=0.38\textwidth]{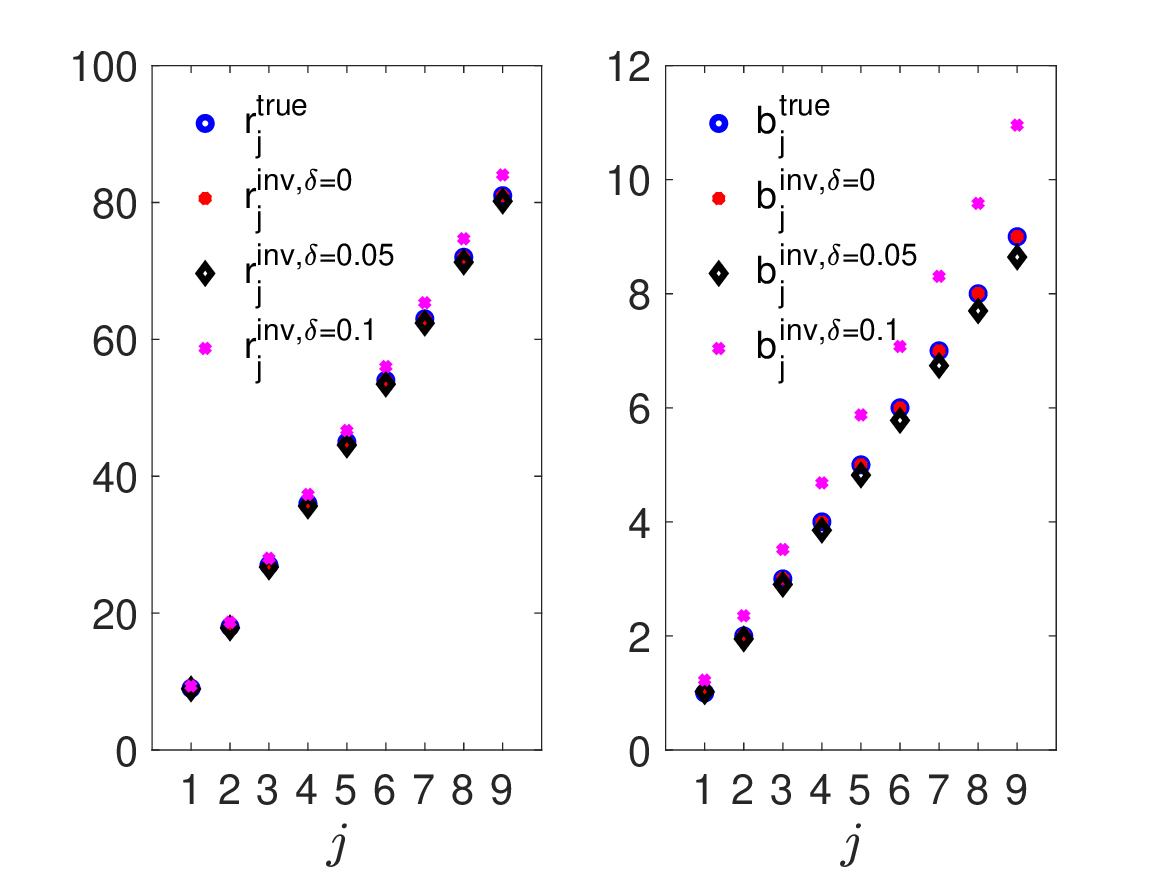}}
&\adjustbox{trim=3.7mm 0mm 5.4mm 2.5mm, clip}{\includegraphics[width=0.38\textwidth]{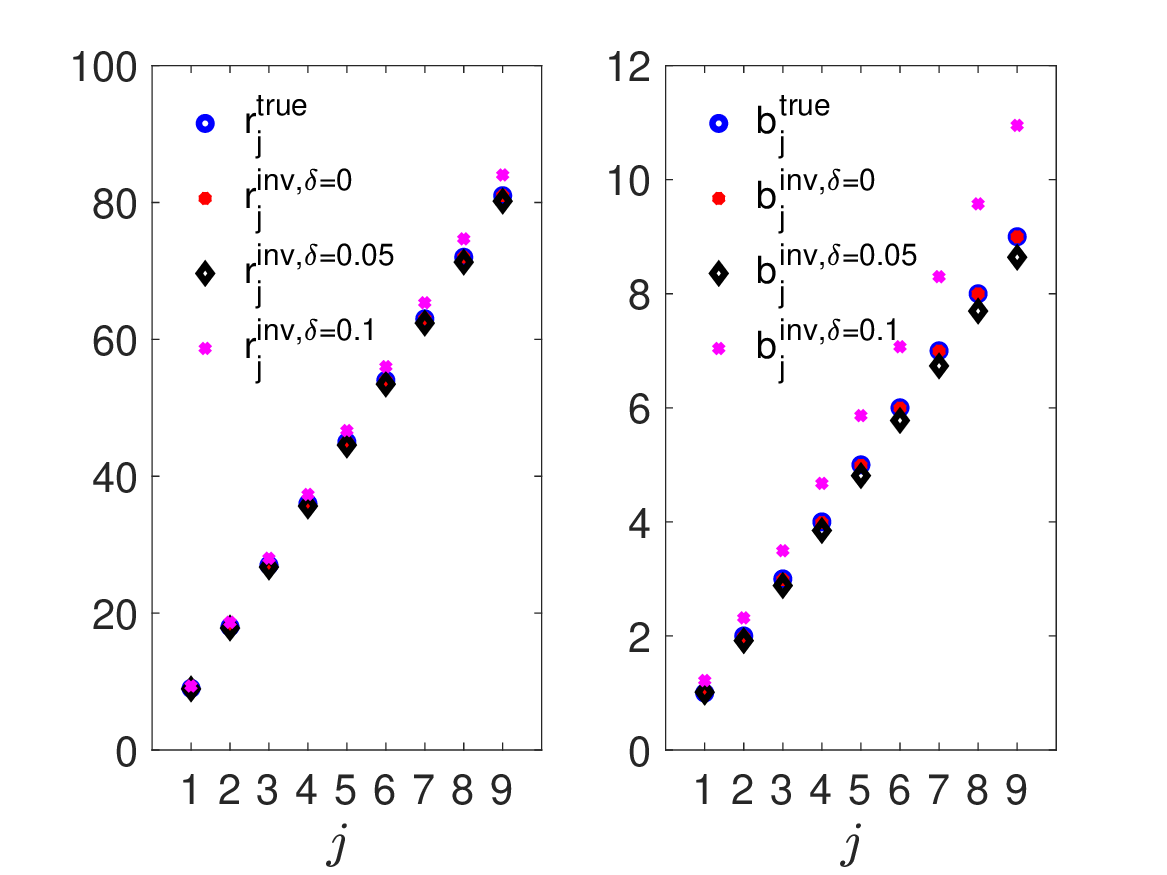}}
\end{tabular}
\caption{Comparisons between reconstructions $(r_j^{inv},\,b_j^{inv},\,D^{inv})$ and true values $(r_j^{true},\,b_j^{true},\,D^{true})$ for $N=9$ and noise level $\delta=0,\,0.05,\,0.1$. Here, each row represents the reconstructions for $D=0.5,\,1,\,5$, respectively. Each column represents the reconstructions for the pair of frequencies $(k_1,\,k_2)=(81,\,91),\,(81,\,501),\,(81,\,1001)$
}
\label{fig_N9}
\end{figure}

\begin{example}\label{exm_2}
Let $N=9$ and $D=0.5,\,1,\,5$. The reconstructions of $(r_j^{inv},\,b_j^{inv}),\,j=1,\,2,\,\cdots,\,N$ and $D^{inv}$ are plotted in Figure \ref{fig_N9} and listed in Table \ref{tbl:N9}, respectively. 
\end{example}

From the results shown in Figure \ref{fig_N9} and Table \ref{tbl:N9}, we can obtain the results similar Example \ref{exm_1} for $N=5$. In other words, the reconstruction scheme is not affected by the dimension.

\section{Conclusions and discussions}
%\noindent

We examined the numerical performance of our proposed inversion scheme. For two given sets of different clustered eigenvalues, the reconstructions are quite accurate for the noise-free data. Furthermore, the reconstructions are robust to the data perturbed by an additional noise. From the results shown in Example \ref{exm_1} and Example \ref{exm_2}, it is worth mentioning that the ill-posedness of reconstruction is not dependent on the dimension. For the glassy state case, $D$ can be also recovered from $r_i$'s and $b_i$'s by \eqref{D} without using \eqref{n3.2}. Even in this way, the reconstructed results are the same for different pairs of frequencies and each noise level, which are equal to the results shown in Table \ref{tbl:N5} and Table \ref{tbl:N9} for $(k_1,\,k_2)=(81,\,1001)$. Hence, the robustness of recovering $D$ against noise follows from those of $r_i$'s and $b_i$'s.

While the number of targets,  $D,\,b_j,\,r_j,j=1,\,\cdot,\, N,$ to recover is $2N+1$, the number of measured data is $2(N+2)$ consisting of the two sets of clustered eigenvlues $\{a_j^k\}_{j=1}^{N+2}$ with $k=1,2$. So, it is natural to ask "Do we really need more measured data to recover less number of targets?". Of course, there is a chance to reduce the number of measured data, but we speculate that the inversion argument becomes mathematically more involved. The biggest advantage of our inversion method is that it is easy to understand even for engineers without knowing any mathematically sophisticated arguments. Also, in practice, it won't be a big deal to use two frequencies and the associated two sets of clustered eigenvalues.

Our next research plan is to collaborate with physicists and engineers to verify the effectiveness of our inversion method on real data and make any necessary improvements.

\subsection*{Acknowledgements}
As for the funding, the second author was supported by the Simons Foundation under the MATH + X program, the National Science Foundation under grant DMS-2108175, and the corporate members of the Geo-Mathematical Imaging Group at Rice University. The third author was supported by NSFC-RGC Joint Research Grant No. 12161160314 and the Natural Science Foundation Innovation Research Team Project of Guangxi (Grant No.2025GXNSFGA069001). The fourth author was partially supported by the Ministry of Science and Technology of Taiwan (Grant No. NSTC 111-2115-M-006-013-MY3). The fifth author was partially supported by JSPS KAKENHI (Grant No. JP22K03366, JP25K07076). 

\subsection*{Data availability}
Data will be made available on request.

\begin{table}[htp]
\centering
\caption{Reconstruction of $D$ from noisy eigenvalues for $N=5$ and different frequencies $(k_1,\,k_2)$}
\label{tbl:N5}
\begin{tabular}{ccccc}
\hline
&$(k_1,\,k_2)$ & $(81,\,91)$   &  $(81,\,501)$ & $(81,\,1001)$   \\ \hline
\multirow{3}*{$D=0.5$}  & $D^{inv,\delta=0}$ &0.5011& 0.5000 & 0.5000\\
                        &$D^{inv,\delta=0.05}$ & 0.4892  & 0.4881 & 0.4881\\
                        &$D^{inv,\delta=0.1}$ & 0.5572& 0.5561& 0.5561\\\hline  
 \multirow{3}*{$D=1$}  & $D^{inv,\delta=0}$ &1.0014& 1.0000 & 1.0000\\
                         &$D^{inv,\delta=0.05}$ & 0.9775  & 0.9762 & 0.9762\\
                         &$D^{inv,\delta=0.1}$ &1.1137  & 1.1122 & 1.1122\\\hline   
\multirow{3}*{$D=5$}  & $D^{inv,\delta=0}$ &5.0026 & 5.0001 & 5.0000\\
                     &$D^{inv,\delta=0.05}$ & 4.8834 & 4.8810 & 4.8809\\
                       &$D^{inv,\delta=0.1}$ & 5.5637  & 5.5610 & 5.5609\\\hline   
\end{tabular}
\end{table}

\begin{table}[htp]
\centering
\caption{Reconstruction of $D$ from noisy eigenvalues for $N=9$ and different frequencies $(k_1,\,k_2)$}
\label{tbl:N9}
\begin{tabular}{ccccc}
\hline
&$(k_1,\,k_2)$ & $(81,\,91)$   &  $(81,\,501)$ & $(81,\,1001)$   \\ \hline
\multirow{3}*{$D=0.5$}  & $D^{inv,\delta=0}$ &0.5032& 0.5001 & 0.5000\\
                        &$D^{inv,\delta=0.05}$ & 0.4913  & 0.4882 & 0.4881\\
                        &$D^{inv,\delta=0.1}$ & 0.5592& 0.5562& 0.5561\\\hline  
 \multirow{3}*{$D=1$}  & $D^{inv,\delta=0}$ &1.0042& 1.0001 & 1.0000\\
                         &$D^{inv,\delta=0.05}$ & 0.9803  & 0.9763 & 0.9762\\
                         &$D^{inv,\delta=0.1}$ &1.1166  & 1.1123 & 1.1122\\\hline   
\multirow{3}*{$D=5$}  & $D^{inv,\delta=0}$ &5.0093 & 5.0003 & 5.0001\\
                     &$D^{inv,\delta=0.05}$ & 4.8901 & 4.8812 & 4.8810\\
                       &$D^{inv,\delta=0.1}$ & 5.5710  & 5.5612 & 5.5610\\\hline  
\end{tabular}
\end{table}


\begin{thebibliography}{99}
%\bibitem{Agr}
%M. Agranovich, Spectral problems in Lipschitz domains, \emph{Journal of Mathematical Sciences}, 190 (2013) pp. 8-33. % https://doi.org/10.1007/s10958-013-1244-6

%\bibitem{CH}
%R. Courant, D. Hilbert, \emph{Methods of Mathematical Physics, vol. 1}, Interscience, New York, 1953.

\bibitem{CDDLN}
S. Chen, M. de Hoop, Y. Deng, C-L. Lin, G. Nakamura, Clustered eigenvalue problem for glassy state relaxation and its inverse problem, arXiv:2509.16714.

\bibitem{de Verdier}
Y. C. de Verdier, Spectre conjoint d'op\'{e}rateurs pseudo-diff\'{e}rentielles qui commutent, Math. Z., 171 (1980) pp. 51-73. %https://doi.org/10.1007/BF01215054

\bibitem{DKLNT}
M. de Hoop, M. Kimura, C-L. Lin, G. Nakamura, K. Tanuma, Anisotropic extended Burgers model, its relaxation tensor and properties of the associated Boltzmann viscoelastic system, SIAM J. Appl. Math.,  85 (2025), 1172-1188.

\bibitem{DLN}
M. de Hoop, C-L Lin, G. Nakamura, Uniform decaying property of solutions for anisotropic viscoelastic systems, Evolution Equations and Control Theory, 16 (2026) pp. 54-79, Doi: 10.3934/eect.2025067.

\bibitem{Guillemin}
V. Guillemin, Some spectral results for the Laplacian on the n-sphere, Adv. Math., 27 (1978) pp. 273-286.

\bibitem{Gurarie}
D. Gurarie, Inverse spectral problem for the 2-sphere Schr\"odinger operators with zonal potentials, Lett. Math. Phys., 16 (1988) pp. 313-323.


%\bibitem{Jackson}
%I. Jackson, Viscoelastic Behaviour from complementary forced-oscillation and microcreep tests, \emph{Minerals}, 9 (2019) 721. % doi. 10.3390/min9120721


\bibitem{LS}
P. Loreti, D. Sforza, Viscoelastic aspects of glass relaxation models, Physica A, 526 (2019) 120768. %https://doi.org/10.1016/j.physa.2019.04.004

\bibitem{MM}
J. Mauro, Y. Mauro, On the prony series representation of the stretched exponetial relaxation, Physica A, 506 (2018) pp. 75-87. %https://doi.org/10.1016/j.physa.2018.04.047

%\bibitem{MZ}
% J. Mauro, E. Zanoto, Two centuries glass research: Historical trends, current status, and grand challenges for the future, \emph{Int. J. Appl. Glass Sci.}, 5 (2014) pp. 313-327.

\bibitem{OR}
J-P. Ortega, F. Rossmannek, Fading memory and the convolution theorem, arXiv:2408.07386v1.

\bibitem{Phillips}
J. Phillips, Stretched exponential relaxation in molecular and electronic glasses, Reports on Progress in Physics, 59 (1996) pp. 1133-1207.


%\bibitem{Jia Shi}
%J. Shi, R. Li, Y. Xi, Y. Saad, M. de Hoop, 
%A Non-perturbative Approach to Computing Seismic Normal Modes in Rotating Planets. \emph{J. Sci. Comput.}, 91 (2022) 67. %https://doi.org/10.1007/s10915-022-01836-5}.


\bibitem{Song}
L. Song, Y. Gao, P. Zou, W. Xu, M. Gao, Y. Zhang, J. Huo, F. Li, J. Qiao, L-M. Wang, J-Q. Wang, Detecting the exponential relaxation spectrum in glasses by high-precision nanocalorimetry, Proc. Natl. Acad. Sci. U.S.A., 120 (2023) e2302776120. %https://doi.org/10.1073/pnas.2302776120.

\bibitem{Weinstein}
A. Weinstein, Asymptotics of eigenvalue clusters for the Laplacian plus a potential, Duke Math. J., 44 (1977) pp. 883-892.

%\bibitem{YP}
%D. Yuen, W. Peltier, Normal modes of the viscoelastic earth, \emph{Geophys. J. R. astr. Soc.}, 69 (1982) pp. 495-526.

%\bibitem{ZYMK}
%T. Zhou, J. Yan, J. Masuda, T. Kuriyagawa, Investigation on the viscoelasticity of optical glass in ultrapresicion lens molding process, J. Materials Processing Tech. 209 (2009) pp.4484-4489.

\end{thebibliography}
\end{document}